# Positivity of Chern classes of Schubert cells and varieties

June Huh

ABSTRACT

We show that the Chern-Schwartz-MacPherson class of a Schubert cell in a Grassmannian is represented by a reduced and irreducible subvariety in each degree. This gives an affirmative answer to a positivity conjecture of Aluffi and Mihalcea.

## 1. Introduction

The classical *Schubert varieties* in the Grassmannian of $d$-planes in a vector space $E$ are among the most studied singular varieties in algebraic geometry. The subject of this paper is the study of *Chern classes* of Schubert cells and varieties.

There is a good theory of Chern classes for singular or noncomplete complex algebraic varieties. If $X^\circ$ is a locally closed subset of a complete variety $X$, then the *Chern-Schwartz-MacPherson class* of $X^\circ$ is an element in the Chow group

$$c_{SM}(X^\circ) \in A_*(X),$$

which agrees with the total homology Chern class of the tangent bundle of $X$ if $X$ is smooth and $X = X^\circ$. The Chern-Schwartz-MacPherson class satisfies good functorial properties which, together with the normalization for smooth and complete varieties, uniquely determines it. Basic properties of the Chern-Schwartz-MacPherson class are recalled in Section 2.1.

If $\underline{\alpha} = (\alpha_1 \geqslant \alpha_2 \geqslant \cdots \geqslant \alpha_d \geqslant 0)$ is a partition, then there is a corresponding Schubert variety $\mathbb{S}(\underline{\alpha})$ in the Grassmannian of $d$-planes in $E$, parametrizing $d$-planes which satisfy incidence conditions with a flag of subspaces determined by $\underline{\alpha}$. See Section 2.2 for our notational conventions. The Schubert variety is a disjoint union of Schubert cells

$$\mathbb{S}(\underline{\alpha}) = \coprod_{\underline{\beta} \leqslant \underline{\alpha}} \mathbb{S}(\underline{\beta})^\circ,$$

where the union is over all $\underline{\beta} = (\beta_1 \geqslant \beta_2 \geqslant \cdots \geqslant \beta_d \geqslant 0)$ which satisfy $\beta_i \leqslant \alpha_i$ for all $i$. Since each Schubert cell $\mathbb{S}(\underline{\beta})^\circ$ is isomorphic to an affine space, the Chow group of $\mathbb{S}(\underline{\alpha})$ is freely generated by the classes of the closures $\big[\mathbb{S}(\underline{\beta})\big]$. Therefore we may write

$$c_{SM}\big(\mathbb{S}(\underline{\alpha})^\circ\big) = \sum_{\underline{\beta} \leqslant \underline{\alpha}} \gamma_{\underline{\alpha},\underline{\beta}} \big[\mathbb{S}(\underline{\beta})\big] \in A_*\big(\mathbb{S}(\underline{\alpha})\big)$$

for uniquely determined coefficients $\gamma_{\underline{\alpha},\underline{\beta}} \in \mathbb{Z}$.

The author was partially supported by NSF grant DMS-0943832.

Various explicit formulas for these coefficients are obtained in [AM09]. One of the formulas says that $\gamma_{\underline{\alpha},\underline{\beta}}$ is the sum of the binomial determinants

$$\gamma_{\underline{\alpha},\underline{\beta}} = \sum_L \det\left[\binom{\alpha_i - l_{i,i+1} - l_{i,i+2} - \cdots - l_{i,d}}{\beta_j + i - j + l_{1,i} + l_{2,i} + \cdots + l_{i-1,i} - l_{i,i+1} - l_{i,i+2} - \cdots - l_{i,d}}\right]_{1 \leqslant i,j \leqslant d}$$

where the sum is over all strictly upper triangular nonnegative integral matrices $L = [l_{p,q}]_{1 \leqslant p < q \leqslant d}$ such that

$$0 \leqslant l_{p,p+1} + l_{p,p+2} + \cdots + l_{p,d} \leqslant \alpha_{p+1} \quad \text{for} \quad 1 \leqslant p < d.$$

For example, $\gamma_{(3 \geqslant 2 \geqslant 1),(2 \geqslant 0 \geqslant 0)}$ is the sum of the determinants of the matrices

$$\begin{pmatrix} 3 & 0 & 0 \\ 0 & 1 & 0 \\ 0 & 1 & 1 \end{pmatrix}, \begin{pmatrix} 2 & 0 & 0 \\ 0 & 2 & 1 \\ 0 & 1 & 1 \end{pmatrix}, \begin{pmatrix} 2 & 0 & 0 \\ 0 & 1 & 0 \\ 0 & 0 & 1 \end{pmatrix}, \begin{pmatrix} 1 & 0 & 0 \\ 0 & 1 & 2 \\ 0 & 1 & 1 \end{pmatrix}, \begin{pmatrix} 1 & 0 & 0 \\ 0 & 2 & 1 \\ 0 & 0 & 1 \end{pmatrix}, \begin{pmatrix} 1 & 0 & 0 \\ 0 & 1 & 0 \\ 0 & 0 & 0 \end{pmatrix},$$

$$\begin{pmatrix} 3 & 0 & 0 \\ 0 & 0 & 0 \\ 0 & 0 & 1 \end{pmatrix}, \begin{pmatrix} 2 & 0 & 0 \\ 0 & 1 & 0 \\ 0 & 0 & 1 \end{pmatrix}, \begin{pmatrix} 2 & 0 & 0 \\ 0 & 0 & 0 \\ 0 & 0 & 0 \end{pmatrix}, \begin{pmatrix} 1 & 0 & 0 \\ 0 & 1 & 1 \\ 0 & 0 & 1 \end{pmatrix}, \begin{pmatrix} 1 & 0 & 0 \\ 0 & 1 & 0 \\ 0 & 0 & 0 \end{pmatrix}, \begin{pmatrix} 1 & 0 & 0 \\ 0 & 0 & 0 \\ 0 & 0 & 0 \end{pmatrix}.$$

That is,

$$\gamma_{(3 \geqslant 2 \geqslant 1),(2 \geqslant 0 \geqslant 0)} = 3 + 2 + 2 + (-1) + 2 + 0 + 0 + 2 + 0 + 1 + 0 + 0 = 11.$$

Based on substantial computer calculations, Aluffi and Mihalcea conjectured that all $\gamma_{\underline{\alpha},\underline{\beta}}$ are nonnegative [AM09, Conjecture 1].

CONJECTURE 1. *For all $\underline{\beta} \leqslant \underline{\alpha}$, the coefficient $\gamma_{\underline{\alpha},\underline{\beta}}$ is nonnegative.*

When $d = 2$, the classical Lindström-Gessel-Viennot lemma shows that $\gamma_{\underline{\alpha},\underline{\beta}}$ is the number of certain nonintersecting lattice paths joining pairs of points in the plane, and hence nonnegative [AM09, Theorem 4.5].

The following is the main result of this paper. Fix a nonnegative integer $k \leqslant \dim \mathbb{S}(\underline{\alpha})$, and write $c_{SM}\big(\mathbb{S}(\underline{\alpha})^\circ\big)_k$ for the $k$-dimensional component of $c_{SM}\big(\mathbb{S}(\underline{\alpha})^\circ\big)$ in $A_k\big(\mathbb{S}(\underline{\alpha})\big)$.

THEOREM 2. *There is a nonempty reduced and irreducible $k$-dimensional subvariety $Z(\underline{\alpha})$ of $\mathbb{S}(\underline{\alpha})$ such that*

$$c_{SM}\big(\mathbb{S}(\underline{\alpha})^\circ\big)_k = \big[Z(\underline{\alpha})\big] \in A_k\big(\mathbb{S}(\underline{\alpha})\big).$$

For details on the subvariety $Z(\underline{\alpha})$, see Section 4. The proof of Theorem 2 is based on an explicit description the Chern class of a vector bundle at the level of cycles. This vector bundle lives on a carefully chosen desingularization of $\mathbb{S}(\underline{\alpha})$, and it is not globally generated in general.

Since any 0-dimensional subvariety is a point, the assertion of Theorem 2 when $k = 0$ is just

$$\chi\big(\mathbb{S}(\underline{\alpha})^\circ\big) = \int_{\mathbb{S}(\underline{\alpha})} c_{SM}\big(\mathbb{S}(\underline{\alpha})^\circ\big) = 1.$$

In general, homology classes representable by a reduced and irreducible subvariety have significantly stronger properties than those representable by an effective cycle. These stronger properties are sometimes of interest in applications [Huh12a, Huh12b]. Unfortunately, little seems to be known about homology classes of subvarieties of a Grassmannian. For the case of curves and multiples of Schubert varieties, however, see [Bry10, Cos11, CR13, Hon05, Hon07, Per02].

It is known that the cone of effective cycles in $A_k\big(\mathbb{S}(\underline{\alpha})\big) \otimes \mathbb{Q}$ is a polyhedral cone generated by the classes of $k$-dimensional $\mathbb{S}(\underline{\beta})$ with $\underline{\beta} \leqslant \underline{\alpha}$ [FMSS95]. Therefore Theorem 2 gives an affirmative answer to Conjecture 1.

COROLLARY 3. *For all $\underline{\beta} \leqslant \underline{\alpha}$, the coefficient $\gamma_{\underline{\alpha},\underline{\beta}}$ is nonnegative.*

Corollary 3 was previously known for all $\underline{\alpha}$ when $d = 2$ [AM09] or $d = 3$ [Mih07], and for all $\underline{\beta} \leqslant \underline{\alpha}$ such that the codimension of $\mathbb{S}(\underline{\beta})$ in $\mathbb{S}(\underline{\alpha})$ is at most 4 [Str11].

It also follows from Theorem 2 that the Chern-Schwartz-MacPherson class of the Schubert variety

$$c_{SM}\big(\mathbb{S}(\underline{\alpha})\big) = \sum_{\underline{\beta} \leqslant \underline{\alpha}} c_{SM}\big(\mathbb{S}(\underline{\beta})^\circ\big)$$

is represented by an effective cycle. This weaker version of positivity was obtained in [Jon10, Theorem 6.5] for a certain infinite class of partitions $\underline{\alpha}$ using Zelevinsky's small resolution.

Finding a positive combinatorial formula for $\gamma_{\underline{\alpha},\underline{\beta}}$ remains as a very interesting problem. As mentioned before, $\gamma_{\underline{\alpha},\underline{\beta}}$ is the number of certain nonintersecting lattice paths joining pairs of points in the plane when $d = 2$. A similar positive combinatorial formula is known for $d = 3$ [Mih07, Corollary 3.10]. The reader will find useful discussions and numerical tables of $\gamma_{\underline{\alpha},\underline{\beta}}$ in [AM09, Mih07, Jon07, Jon10, Str11, Web12].

ACKNOWLEDGEMENTS

The author is grateful to Dave Anderson, William Fulton, and Bernd Sturmfels for useful comments. He thanks Mircea Mustaţă for helpful discussions.

## 2. Preliminaries

**2.1**

We briefly recall the basic properties of the Chern-Schwartz-MacPherson class. More details can be found in [Alu05, Ken90, Mac74, Sch05].

Let $X$ be a complete complex algebraic variety. The group of constructible functions on $X$ is the free abelian group $C(X)$ generated by functions of the form

$$\mathbf{1}_W = \begin{cases} 1, & x \in W, \\ 0, & x \notin W, \end{cases}$$

where $W$ is a closed subvariety of $X$. If $f : X \longrightarrow Y$ is a morphism between complete varieties, then the pushforward $f_*$ is defined to be the homomorphism

$$f_* : C(X) \longrightarrow C(Y), \qquad \mathbf{1}_W \longmapsto \Big(y \longmapsto \chi\big(f^{-1}(y) \cap W\big)\Big)$$

where $\chi$ stands for the topological Euler characteristic. This defines a functor $C$ from the category of complete varieties to the category of abelian groups.

DEFINITION 4. *The Chern-Schwartz-MacPherson class is the unique natural transformation*

$$c_{SM} : C \longrightarrow A_*$$

*such that*

$$c_{SM}(\mathbf{1}_X) = c(T_X) \cap [X] \in A_*(X)$$

*if $X$ is a smooth and complete variety with the tangent bundle $T_X$.* When $X^\circ$ is a locally closed subset of $X$, we write

$$c_{SM}(X^\circ) := c_{SM}(\mathbf{1}_{X^\circ}).$$

The functoriality of $c_{SM}$ says that, for any $f : X \longrightarrow Y$ as above, we have the commutative diagram

$$\begin{array}{ccc} C(X) & \xrightarrow{c_{SM}} & A_*(X) \\ {\scriptstyle f_*}\downarrow & & \downarrow{\scriptstyle f_*} \\ C(Y) & \xrightarrow[c_{SM}]{} & A_*(Y). \end{array}$$

The uniqueness of $c_{SM}$ follows from the functoriality, the resolution of singularities, and the normalization for smooth and complete varieties. The existence of $c_{SM}$, which was once a conjecture of Deligne and Grothendieck, was proved by MacPherson in [Mac74]. The Chern-Schwartz-MacPherson class satisfies the inclusion-exclusion formula

$$c_{SM}(\mathbf{1}_{U_1 \cup U_2}) = c_{SM}(\mathbf{1}_{U_1}) + c_{SM}(\mathbf{1}_{U_2}) - c_{SM}(\mathbf{1}_{U_1 \cap U_2})$$

and captures the topological Euler characteristic as its degree

$$\chi(U) = \int c_{SM}(\mathbf{1}_U).$$

Here $U, U_1, U_2$ can be any constructible subset of a complete variety. For a construction of $c_{SM}$ with an emphasis on noncomplete varieties, see [Alu06a, Alu06b].

**2.2**

We define the Schubert variety $\mathbb{S}(\underline{\alpha})$ corresponding to a partition $\underline{\alpha}$ in the Grassmannian of $d$-planes $\mathrm{Gr}_d(E)$. Schubert varieties will only appear at the last section of this paper.

Our notation for Schubert varieties is consistent with that of [AM09]. In the study of homology Chern classes, this 'homological' notation has advantages over the more common 'cohomological' notation.

Let $E$ be a complex vector space with an ordered basis $e_1, \ldots, e_{n+d}$, and take $F_k$ to be the subspace spanned by the first $k$ vectors in this basis.

DEFINITION 5. Let $\underline{\alpha} = (\alpha_1 \geqslant \alpha_2 \geqslant \cdots \geqslant \alpha_d \geqslant 0)$ be a partition with $n \geqslant \alpha_1$.

(i) The *Schubert variety* corresponding to $\underline{\alpha}$ is the subvariety

$$\mathbb{S}(\underline{\alpha}) := \Big\{ V \mid \dim(V \cap F_{\alpha_{d+1-i}+i}) \geqslant i \text{ for } i = 1, \ldots, d \Big\} \subseteq \mathrm{Gr}_d(E).$$

(ii) The *Schubert cell* corresponding to $\underline{\alpha}$ is the open subset of $\mathbb{S}(\underline{\alpha})$

$$\mathbb{S}(\underline{\alpha})^\circ := \Big\{ V \mid \dim(V \cap F_{\alpha_{d+1-i}+i}) = i, \ \dim(V \cap F_{\alpha_{d+1-i}+i-1}) = i-1 \text{ for } i = 1, \ldots, d \Big\}.$$

We summarize the main properties of Schubert cells and varieties:

1. Writing $\underline{\beta} \leqslant \underline{\alpha}$ for the ordering $\beta_i \leqslant \alpha_i$ for all $i$, we have

$$\mathbb{S}(\underline{\alpha})^\circ = \mathbb{S}(\underline{\alpha}) \setminus \bigcup_{\underline{\beta} < \underline{\alpha}} \mathbb{S}(\underline{\beta}).$$

2. The Schubert cell $\mathbb{S}(\underline{\alpha})^\circ$ is isomorphic to the affine space $\mathbb{C}^{\alpha_1 + \cdots + \alpha_d}$.

3. The Schubert cell $\mathbb{S}(\underline{\alpha})^\circ$ is an orbit under the natural action of $B$ on $\mathrm{Gr}_d(E)$.

Here $B$ is the subgroup of the general linear group of $E$ which consists of all invertible upper triangular matrices with respect to the ordered basis $e_1, \ldots, e_{n+d}$. The reader will find details in [AM09, Bri05, Ful97].

## 3. Chern classes of almost homogeneous varieties

In this section, $B$ is a connected affine algebraic group with the Lie algebra $\mathfrak{b}$.

### 3.1

Suppose $B$ acts on an irreducible projective variety $Y$ with an open dense orbit $Y^\circ$. We say that $Y$ is *almost homogeneous* with respect to the action of $B$. For example, $Y$ can be the Schubert variety $\mathbb{S}(\underline{\alpha})$ of the previous section.

DEFINITION 6. A *$B$-finite log-resolution* of $Y$ is a proper $B$-equivariant map $\pi : X \longrightarrow Y$ such that

(i)  $X$ is smooth and has finitely many $B$-orbits,
(ii) $\pi^{-1}(Y^\circ) \longrightarrow Y^\circ$ is an isomorphism, and
(iii) the complement of $\pi^{-1}(Y^\circ)$ in $X$ is a divisor with normal crossings.

The main result of this section is the following sufficient condition for the Chern-Schwartz-MacPherson class of an almost homogeneous $B$-variety to be effective.

THEOREM 7. *Suppose $Y$ has a $B$-finite log-resolution. Then there are subvarieties $Z_1, \ldots, Z_p$ of $Y$ and nonnegative integers $n_1, \ldots, n_p$ such that*

$$c_{SM}(Y^\circ) = \sum_{i=1}^{p} n_i [Z_i] \in A_*(Y).$$

In short, the Chern-Schwartz-MacPherson class of $Y^\circ$ is represented by an effective cycle on $Y$ if $Y$ has a $B$-finite resolution.

When $Y$ is the Schubert variety $\mathbb{S}(\underline{\alpha})$, the conclusion of Theorem 7 is much weaker than that of Theorem 2. However, the main construction which leads to the proof of Theorem 7 will be essential in the proof of Theorem 2.

The rest of this section is devoted to the proof of Theorem 7.

### 3.2

As a preparation, we recall basic results on algebraic groups actions and algebraic vector fields. General references are [MO67] and [Ram64].

Suppose $B$ acts on a smooth and irreducible projective variety $X$. There is an algebraic group homomorphism from $B$ to the connected automorphism group

$$L : B \longrightarrow \mathrm{Aut}^\circ(X), \qquad b \longmapsto \big(x \longmapsto b \cdot x\big).$$

The differential of $L$ at the identity is the *Lie homomorphism* between the Lie algebras

$$\mathfrak{b} \longrightarrow \Gamma(X, T_X).$$

Explicitly, the Lie homomorphism maps $\xi \in \mathfrak{b}$ to the corresponding fundamental vector field

$$x \longmapsto \frac{d}{dt}\bigg|_{t=0} \Big( \exp(-t\xi) \cdot x \Big).$$

If we define the $B$-action on the vector fields on $X$ by

$$\big(x \longmapsto v(x)\big) \longmapsto \big(x \longmapsto d(b \cdot -) v(b^{-1} \cdot x)\big),$$

then the Lie homomorphism is $B$-equivariant with respect to the adjoint action of $B$ on $\mathfrak{b}$. Evaluating the Lie homomorphism, we have the homomorphism between the $B$-linearized vector bundles
$$\mathscr{L}_X : \mathfrak{b}_X \longrightarrow T_X,$$
where $\mathfrak{b}_X$ is the trivial vector bundle on $X$ modeled on $\mathfrak{b}$.

### 3.3

Let $S$ be an orbit of the $B$-action on $X$, and write $\iota$ for the inclusion $S \longrightarrow X$. A choice of a base point $x_0 \in S$ defines the orbit map
$$B \longrightarrow S, \qquad b \longmapsto b \cdot x_0.$$
This identifies $S$ with $B/H$, where $H$ is the isotropy group $B_{x_0}$. The Lie homomorphism
$$\mathfrak{b} \longrightarrow \Gamma(S, T_S)$$
gives the $B$-linearized vector bundle homomorphism
$$\mathscr{L}_S : \mathfrak{b}_S \longrightarrow T_S,$$
and $\mathscr{L}_S$ fits into the commutative diagram

$$\begin{array}{ccc} \mathfrak{b}_S & \xrightarrow{\mathscr{L}_S} & T_S \\ {\scriptstyle \mathscr{L}_X|_S}\downarrow & \swarrow {\scriptstyle \iota_*} & \\ T_X|_S. & & \end{array}$$

Over the base point $x_0$, $\mathscr{L}_S$ can be identified with the surjective linear map
$$\mathfrak{b} \longrightarrow \mathfrak{b}/\mathfrak{h},$$
where $\mathfrak{h}$ is the Lie algebra of $H$. Since $S$ is homogeneous, $\mathscr{L}_S$ is surjective over every point of $S$, and $\ker(\mathscr{L}_S)$ is a vector bundle over $S$.

DEFINITION 8. The *bundle of isotropy Lie algebras* over $S$ is the locally closed subset
$$\Sigma_S := \mathbb{P}\big(\ker(\mathscr{L}_S)\big) \subseteq X \times \mathbb{P}(\mathfrak{b}).$$

Note that $\Sigma_S$ is a smooth and irreducible closed subset of $S \times \mathbb{P}(\mathfrak{b})$. We denote the two projections by

$$\begin{array}{ccc} & \Sigma_S & \\ {\scriptstyle \mathrm{pr}_{1,S}}\swarrow & & \searrow {\scriptstyle \mathrm{pr}_{2,S}} \\ S & & \mathbb{P}(\mathfrak{b}). \end{array}$$

If we write $\mathfrak{b}_x$ for the Lie algebra of the isotropy group $B_x$, then
$$\Sigma_S = \Big\{ (x, \xi) \mid x \in S \text{ and } \xi \in \mathfrak{b}_x \Big\}.$$
The dimension of $\Sigma_S$ is equal to the dimension of $\mathbb{P}(\mathfrak{b})$, independent of the dimension of $S$.

### 3.4

Let $D$ be a simple normal crossings divisor on $X$. The *logarithmic tangent sheaf* of $(X, D)$ is the subsheaf of the tangent sheaf
$$\mathcal{T}_X(-\log D) \subseteq \mathcal{T}_X$$

consisting of those derivations which preserve the ideal sheaf $\mathcal{O}_X(-D)$. Since $D$ is a divisor with simple normal crossings, the logarithmic tangent sheaf is locally free of rank equal to the dimension of $X$. General references on logarithmic tangent sheaves are [Del70] and [Sai80].

We write $T_X(-\log D)$ for the logarithmic tangent bundle, the vector bundle corresponding to the logarithmic tangent sheaf. The following equality follows from a construction of the Chern-Schwartz-MacPherson class [Alu06a, Alu06b].

THEOREM 9. *We have*

$$c_{SM}(\mathbf{1}_{X \setminus D}) = c\bigl(T_X(-\log D)\bigr) \cap [X] \in A_*(X).$$

For precursors, see [Alu99, GP02] and also Schwartz's construction of the Chern-Schwartz-MacPherson class [BSS09, Sch65a, Sch65b]. Our goal is to show that $X$ has enough logarithmic vector fields so as to make the right-hand side of Theorem 9 effective when $D$ is $B$-invariant and $X$ has finitely many $B$-orbits.

Suppose from now on that $D$ is invariant under the action of $B$. This implies that the Lie homomorphism of Section 3.2 factors through

$$\mathscr{L} : \mathfrak{b} \longrightarrow \Gamma\bigl(X, T_X(-\log D)\bigr).$$

Evaluating the sections, we have the homomorphism between $B$-linearized vector bundles

$$\mathscr{L}_{X,D} : \mathfrak{b}_X \longrightarrow T_X(-\log D).$$

We denote the induced linear map between the fibers over $x \in X$ by

$$\mathscr{L}_{X,D,x} : \mathfrak{b} \longrightarrow T_{X,x}(-\log D).$$

DEFINITION 10. The *variety of critical points* of $(X, D)$ is the closed subset

$$\mathfrak{X} := \Bigl\{(x, \xi) \mid \mathscr{L}_{X,D,x}(\xi) = 0\Bigr\} \subseteq X \times \mathbb{P}(\mathfrak{b}).$$

We denote the two projections by

$$\begin{array}{ccc} & \mathfrak{X} & \\ {}^{\mathrm{pr}_1}\swarrow & & \searrow {}^{\mathrm{pr}_2} \\ X & & \mathbb{P}(\mathfrak{b}). \end{array}$$

The first projection $\mathrm{pr}_1 : \mathfrak{X} \longrightarrow X$ may not be a projective bundle, but the restriction $\mathrm{pr}_1^{-1}(S) \longrightarrow S$ is a projective bundle for each $B$-orbit $S$ in $X$. These projective bundles have different ranks in general.

*Remark* 11. When $\mathscr{L}_{X,D}$ is surjective, the pair $(X, D)$ is said to be *log-homogeneous* under the action of $B$ [Bri07]. In this case, $\mathfrak{X}$ is the projectivization of the vector bundle denoted by $R_X$ in [Bri09, Section 2].

For log-homogeneous varieties, the conclusion of Theorem 7 is a standard fact [Ful98, Example 12.1.7]. However, in our main case of interest, $(X, D)$ is rarely log-homogeneous under $B$. In fact, if $(X, D)$ is log-homogeneous under a *solvable* affine algebraic group $B$, then $X$ should be a toric variety of a maximal torus $T \subseteq B$ [Bri07, Theorem 3.2.1].

We refer to [BJ08, BK05, Kir06, Kir07] for studies of Chern classes of the logarithmic tangent bundle of log-homogeneous varieties.

**3.5**

Define $X_0 := X$, $X_1 := D$, and a sequence of closed subsets

$$X_0 \supsetneq X_1 \supsetneq X_2 \supsetneq X_3 \supsetneq \cdots \quad \text{where} \quad X_{i+1} := \text{Sing}(X_i) \quad \text{for} \quad i \geqslant 1.$$

We introduce two decompositions of $X$ into smooth locally closed subsets, the orbit decomposition $\mathcal{S}_{\text{orb}}$ and the singular decomposition $\mathcal{S}_{\text{sing}}$:

$$\mathcal{S}_{\text{orb}} := \big\{ S \mid S \text{ is a } B\text{-orbit in } X \big\},$$
$$\mathcal{S}_{\text{sing}} := \big\{ S \mid S \text{ is a connected component of some } X_i \setminus X_{i+1} \big\}.$$

Since $B$ is connected and $D$ is invariant under the action of $B$, the orbit decomposition refines the singular decomposition. We write the variety of critical points as a disjoint union by taking inverse images over the $B$-orbits in $X$:

$$\mathfrak{X} = \coprod_{S \in \mathcal{S}_{\text{orb}}} \mathfrak{X}_S \quad \text{where} \quad \mathfrak{X}_S := \text{pr}_1^{-1}(S).$$

As in Section 3.3, we denote the bundle of isotropy Lie algebras over $S$ by $\Sigma_S$.

LEMMA 12. *$\mathfrak{X}_S$ is a closed subset of $\Sigma_S$ for each $B$-orbit $S$ in $X$.*

*Proof.* Let $S'$ be the unique element of $\mathcal{S}_{\text{sing}}$ containing $S$. Any section of $T_X(-\log D)$ preserves the ideal sheaf of $S'$ and defines a derivation of $\mathcal{O}_{S'}$. Denote the corresponding vector bundle homomorphism over $S'$ by

$$\varphi : T_X(-\log D)|_{S'} \longrightarrow T_{S'}.$$

Note that the restriction of $\varphi$ to $S$ fits into the commutative diagram:

$$\begin{array}{ccc}
 & & 0 \\
 & & \downarrow \\
\mathfrak{b}_S & \xrightarrow{\mathscr{L}_S} & T_S \\
{\scriptstyle \mathscr{L}_{X,D}|_S} \downarrow & & \downarrow {\scriptstyle \iota_*} \\
T_X(-\log D)|_S & \xrightarrow[\varphi|_S]{} & T_{S'}|_S.
\end{array}$$

Here $\mathscr{L}_S$ is the vector bundle homomorphism of Section 3.3, $\mathscr{L}_{X,D}|_S$ is the restriction to $S$ of the vector bundle homomorphism of Section 3.4, and $\iota_*$ is the differential of the inclusion $\iota : S \to S'$. Since $\iota_*$ is injective, $\mathscr{L}_{X,D,x}(\xi) = 0$ implies $\mathscr{L}_{S,x}(\xi) = 0$ for any $x \in S$ and $\xi \in \mathfrak{b}$. □

**3.6**

*Proof of Theorem 7.* Choose a $B$-finite log-resolution $\pi : X \longrightarrow Y$ and define $X^\circ := \pi^{-1}(Y^\circ)$. By the functoriality of the Chern-Schwartz-MacPherson class, we have

$$\pi_* c_{SM}(X^\circ) = c_{SM}(Y^\circ) \in A_*(Y).$$

Since any effective cycle push-forwards to an effective cycle, it is enough to prove that $c_{SM}(X^\circ)$ is represented by an effective cycle on $X$.

Let $D$ be the boundary divisor $X \setminus X^\circ$, and let $k$ be a nonnegative integer less than $\dim X$. Our aim is to show that the $k$-th Chern class

$$c_{SM}(X^\circ)_k = c_{\dim X - k}\big(T_X(-\log D)\big) \cap [X] \in A_k(X)$$

is represented by an effective $k$-cycle.

We recall from Section 3.4 the variety of critical points $\mathfrak{X}$ and the two projections

$$\begin{array}{c} & \mathfrak{X} & \\ {}_{\mathrm{pr}_1}\swarrow & & \searrow{}^{\mathrm{pr}_2} \\ X & & \mathbb{P}(\mathfrak{b}). \end{array}$$

By Lemma 12, we have

$$\mathfrak{X} = \coprod_{S \in \mathcal{S}_{\mathrm{orb}}} \mathfrak{X}_S \subseteq \coprod_{S \in \mathcal{S}_{\mathrm{orb}}} \Sigma_S.$$

Note that each $\Sigma_S$ is irreducible of dimension equal to that of $\mathbb{P}(\mathfrak{b})$. Since $X$ has finitely many $B$-orbits, this shows that each irreducible component of $\mathfrak{X}$ has dimension at most $\dim \mathbb{P}(\mathfrak{b})$.

Let $\Lambda$ be a $(k+1)$-dimensional subspace of $\mathfrak{b}$. If $\Lambda$ is spanned by $\xi_0, \ldots, \xi_k$, then the $k$-th Chern class of $T_X(-\log D)$ is represented by a cycle supported on the locus

$$\mathfrak{D}_k(\Lambda) := \left\{ x \in X \mid \mathscr{L}(\xi_0), \ldots, \mathscr{L}(\xi_k) \text{ are linearly dependent at } x \right\},$$

where $\mathscr{L} : \mathfrak{b} \longrightarrow \Gamma\bigl(X, T_X(-\log D)\bigr)$ is the Lie homomorphism. See [Ful98, Chapter 14]. As a scheme, $\mathfrak{D}_k(\Lambda)$ is defined by $(k+1)$-minors of the matrices for the vector bundle homomorphism

$$\Lambda_X \longrightarrow T_X(-\log D)$$

obtained by restricting $\mathscr{L}_{X,D}$. Set-theoretically,

$$\mathfrak{D}_k(\Lambda) = \mathrm{pr}_1\Bigl(\mathrm{pr}_2^{-1}\bigl(\mathbb{P}(\Lambda)\bigr)\Bigr).$$

We recall the following facts on degeneracy loci from [Ful98, Theorem 14.4]:

(i) Each irreducible component of $\mathfrak{D}_k(\Lambda)$ has dimension at least $k$.
(ii) If all the irreducible components of $\mathfrak{D}_k(\Lambda)$ have dimension $k$, then the Chern class

$$c_{\dim X - k}\bigl(T_X(-\log D)\bigr) \cap [X] \in A_k(X)$$

is represented by a positive cycle supported on $\mathfrak{D}_k(\Lambda)$.

Therefore it is enough to show that all the irreducible components of $\mathfrak{D}_k(\Lambda)$ have dimension at most $k$ for a suitable choice of $\Lambda$.

In fact, all the irreducible components of $\mathrm{pr}_2^{-1}\bigl(\mathbb{P}(\Lambda)\bigr)$ have dimension at most $k$ for a sufficiently general choice of $\Lambda$. This is a general fact on maps of the form

$$\mathfrak{X} \longrightarrow \mathbb{P}^n,$$

where all the irreducible components of $\mathfrak{X}$ has dimension $\leqslant n$. One may argue by induction on $n$, where in the induction step one chooses a hyperplane of $\mathbb{P}^n$ which does not contain the image of any irreducible component of $\mathfrak{X}$. $\square$

Since each irreducible component of the degeneracy locus $\mathfrak{D}_k(\Lambda)$ has dimension at least $k$, the above argument shows that each component of $\mathfrak{D}_k(\Lambda)$ has dimension exactly $k$ for a sufficiently general $\Lambda$. Each of these components is projected from an irreducible component of $\mathfrak{X}$ of maximum possible dimension, namely the dimension of $\mathbb{P}(\mathfrak{b})$. For a later use, we record here this refined conclusion of our analysis.

COROLLARY 13. *The following hold for a sufficiently general $(k+1)$-dimensional subspace $\Lambda \subseteq \mathfrak{b}$.*

(i) Each irreducible component of $\mathfrak{D}_k(\Lambda)$ has the expected dimension $k$.
(ii) Each irreducible component of $\mathfrak{D}_k(\Lambda)$ is the closure of a subvariety of a $B$-orbit $S$ such that $\mathfrak{X}_S = \Sigma_S$.

We express (ii) by saying that the irreducible component of $\mathfrak{D}_k(\Lambda)$ is generically supported on $S$.

There is at least one orbit with the equality $\mathfrak{X}_S = \Sigma_S$, the open dense orbit $S = X^\circ$. Any irreducible component of $\mathfrak{D}_k(\Lambda)$ generically supported on $X^\circ$ will be called *standard*. All the other irreducible components are *exceptional*.

## 4. Irreducibility

In this section, we specialize to the case when $B$ is a Borel subgroup of a connected reductive group $G$. We make use of the following consequence of the strengthened assumption:

– the centralizer of a maximal torus in $B$ is the maximal torus.

Since the union of Cartan subgroups of $B$ contains an open dense subset, it follows that

– the set of semisimple elements of $B$ contains an open dense subset of $B$, and
– the set of semisimple elements of $\mathfrak{b}$ contains an open dense subset of $\mathfrak{b}$.

We will use [Bor91] as a general reference. For Cartan subgroups and Cartan subalgebras, see [TY05, Chapter 29].

Let $P$ be a parabolic subgroup of $G$ containing $B$, and let $Y$ be the closure of a $B$-orbit $Y^\circ$ in $G/P$.

### 4.1

An element $\xi \in \mathfrak{b}$ is said to be *regular* if its centralizer is a Cartan subalgebra of $\mathfrak{b}$. The set of regular elements is open and dense in $\mathfrak{b}$.

DEFINITION 14. A *regular log-resolution* of $Y$ is a proper map $\pi : X \longrightarrow Y$ such that

(i) $\pi : X \longrightarrow Y$ is a $B$-finite log-resolution of $Y$, and
(ii) the isotropy Lie algebra $\mathfrak{b}_x$ contains a regular element of $\mathfrak{b}$ for each $x \in X$.

Of course, it is enough to require the second condition for any one point from each $B$-orbit of $X$.

The following is the main result of this section. Fix a nonnegative integer $k \leqslant \dim Y$, and write $c_{SM}(Y^\circ)_k$ for the $k$-dimensional component of $c_{SM}(Y^\circ)$.

THEOREM 15. *Suppose $Y$ has a regular log-resolution. Then there is a nonempty reduced and irreducible $k$-dimensional subvariety $Z$ of $Y$ such that*

$$c_{SM}(\mathbf{1}_{Y^\circ})_k = [Z] \in A_k(Y).$$

*The subvariety $Z$ can be chosen to be the closure in $Y$ of the locus*

$$Z^\circ(\Lambda) = \Big\{ y \in Y^\circ \mid \Lambda \cap \mathfrak{b}_y \neq 0 \Big\},$$

*where $\Lambda$ is a sufficiently general $(k+1)$-dimensional subspace of $\mathfrak{b}$.*

We will see in Section 5 that the classical Schubert variety $\mathbb{S}(\underline{\alpha})$ has a regular log-resolution. The rest of this section is devoted to the proof of Theorem 15.

**4.2**

Let $S$ be a homogeneous $B$-space. Recall from Section 3.3 the bundle of isotropy Lie algebras

$$\Sigma_S = \left\{(x, \xi) \mid \xi \in \mathfrak{b}_x\right\} \subseteq S \times \mathbb{P}(\mathfrak{b}).$$

We choose a base point $x_0$ and identify $S$ with $B/H$, where $H$ is the isotropy group $B_{x_0}$ with the Lie algebra $\mathfrak{h}$. The *rank* of an affine algebraic group is the dimension of a maximal torus.

LEMMA 16. *If $\mathrm{rank}(B) = \mathrm{rank}(H)$, then*

$$\mathrm{pr}_{2,S} : \Sigma_S \longrightarrow \mathbb{P}(\mathfrak{b}), \qquad (x, \xi) \longmapsto \xi$$

*is a dominant morphism.*

*Proof.* The set of semisimple elements in $\mathfrak{b}$ contains an open dense subset of $\mathfrak{b}$ in our setting. We find a point in $\Sigma_S$ which maps to the class of a given nonzero semisimple element $\xi$ in $\mathbb{P}(\mathfrak{b})$.

Since $\xi$ is semisimple, $\xi$ is tangent to a torus [Bor91, Proposition 11.8]. We may assume that this torus $T_1$ is a maximal torus of $B$.

Let $T_2$ be a maximal torus of $H$. Then $T_2$ is a maximal torus of $B$ because $\mathrm{rank}(B) = \mathrm{rank}(H)$. Since any two maximal tori of $B$ are conjugate, there is an element $b \in B$ such that $T_1 = bT_2b^{-1}$. We have

$$\xi \in \mathfrak{t}_1 = \mathrm{Ad}(b) \cdot \mathfrak{t}_2 \subseteq \mathrm{Ad}(b) \cdot \mathfrak{h} = \mathfrak{b}_{b \cdot x_0}.$$

Therefore $b \cdot x_0$ gives a point in the fiber of $\xi$. $\square$

**4.3**

*Remark* 17. The results of this subsection are not needed for the proof of Theorem 15 if $Y$ is the classical Schubert variety $\mathbb{S}(\underline{\alpha})$.

Let $\Lambda$ be a $(k+1)$-dimensional subspace of $\mathfrak{b}$, and let $\Lambda_r$ be the set of regular elements of $\mathfrak{b}$ in $\Lambda$. Define

$$D_k(\Lambda) := \{x \in S \mid \Lambda \cap \mathfrak{b}_x \neq 0\} \quad \text{and} \quad D_k(\Lambda_r) := \{x \in S \mid \Lambda_r \cap \mathfrak{b}_x \neq 0\}.$$

In terms of the diagram

$$\begin{array}{ccc} & \Sigma_S & \\ \mathrm{pr}_{1,S} \swarrow & & \searrow \mathrm{pr}_{2,S} \\ S & & \mathbb{P}(\mathfrak{b}), \end{array}$$

we have

$$D_k(\Lambda) = \mathrm{pr}_{1,S}\left(\mathrm{pr}_{2,S}^{-1}\big(\mathbb{P}(\Lambda)\big)\right) \quad \text{and} \quad D_k(\Lambda_r) = \mathrm{pr}_{1,S}\left(\mathrm{pr}_{2,S}^{-1}\big(\mathbb{P}(\Lambda_r)\big)\right).$$

Since $\dim \Sigma_S = \dim \mathbb{P}(\mathfrak{b})$, $D_k(\Lambda)$ is either empty or of pure dimension $k$ for a sufficiently general $\Lambda$.

LEMMA 18. *Suppose $\mathfrak{h}$ contains a regular element of $\mathfrak{b}$. Then $D_k(\Lambda_r)$ contains an open dense subset of $D_k(\Lambda)$ for a sufficiently general $\Lambda \subseteq \mathfrak{b}$.*

*Proof.* Note that
$$\mathrm{pr}_{2,S}(\Sigma_S) = \bigcup_{x \in S} \mathbb{P}(\mathfrak{b}_x).$$
The closure of this set is an irreducible subvariety of $\mathbb{P}(\mathfrak{b})$, say $V$. Let $U \subseteq V$ be the open subset of (the classes of) regular elements in $V$. This set $U$ is nonempty by our assumption on $\mathfrak{h}$, and hence $U$ is dense in $V$.

(i) $\dim V \leqslant \mathrm{codim}(\Lambda \subseteq \mathfrak{b})$: In this case, for a sufficiently general $\Lambda$,
$$V \cap \mathbb{P}(\Lambda) = U \cap \mathbb{P}(\Lambda).$$
Therefore $\mathrm{pr}_{2,S}^{-1}(U \cap \mathbb{P}(\Lambda)) = \mathrm{pr}_{2,S}^{-1}(\mathbb{P}(\Lambda))$.

(ii) $\dim V > \mathrm{codim}(\Lambda \subseteq \mathfrak{b})$: In this case, $\mathrm{pr}_{2,S}^{-1}(\mathbb{P}(\Lambda))$ is irreducible for a sufficiently general $\Lambda$ by Bertini's theorem [Laz04, Theorem 3.3.1]. Therefore $\mathrm{pr}_{2,S}^{-1}(U \cap \mathbb{P}(\Lambda))$ is open and dense in $\mathrm{pr}_{2,S}^{-1}(\mathbb{P}(\Lambda))$.

In either case, we see that $D_k(\Lambda_r)$ contains an open dense subset of $D_k(\Lambda)$. $\square$

Let $p$ be a $B$-equivariant morphism between homogeneous $B$-spaces
$$p : S \simeq B/H \longrightarrow B/K, \qquad H \subseteq K \subseteq B.$$
The following lemma can be found in [Kir07, Lemma 3.1].

LEMMA 19. *If $\mathfrak{h}$ contains a regular element of $\mathfrak{b}$ and $\mathrm{rank}(H) < \mathrm{rank}(K)$, then*
$$\dim D_k(\Lambda) > \dim p(D_k(\Lambda))$$
*for a sufficiently general $\Lambda \subseteq \mathfrak{b}$.*

*Proof.* By Lemma 18, $D_k(\Lambda_r)$ contains an open dense subset $D^\circ$ of $D_k(\Lambda)$. It is enough to show that
$$\dim\left(D_k(\Lambda) \cap p^{-1}(p(x))\right) > 0 \quad \text{for all } x \in D^\circ.$$
Let $x$ be a point in $D^\circ$. Since regular elements are semisimple in our setting, there is a nonzero semisimple element $\xi$ in $\Lambda \cap \mathfrak{b}_x \subseteq \mathfrak{b}_{p(x)}$. Choose a maximal torus $T$ of $B_{p(x)}$ tangent to $\xi$ [Bor91, Proposition 11.8].

The maximal torus $T$ is contained in the centralizer of $\xi$ because global and infinitesimal centralizers correspond [Bor91, Section 9.1]. Therefore, for any $t \in T$,
$$\xi = \mathrm{Ad}(t) \cdot \xi \in \Lambda \cap \mathfrak{b}_{t \cdot x} \neq 0.$$
This shows that
$$T \cdot x \subseteq D_k(\Lambda).$$
Since $T$ is contained in $B_{p(x)}$, we have
$$T \cdot x \subseteq D_k(\Lambda) \cap p(p^{-1}(x)).$$
We check that $T \cdot x$ has a positive dimension. If otherwise, $T \cdot x = x$ because $T \cdot x$ is connected. Therefore $T \subseteq B_x$, and this contradicts the assumption that $\mathrm{rank}(H) < \mathrm{rank}(K)$. $\square$

**4.4**

We begin the proof of Theorem 15. Choose a regular log-resolution $\pi : X \longrightarrow Y$ and set

$$X^\circ := \pi^{-1}(Y^\circ), \qquad D := X \setminus X^\circ.$$

By the functoriality, we have

$$\pi_* c_{SM}(X^\circ) = c_{SM}(Y^\circ) \in A_*(Y).$$

Let $\Lambda \subseteq \mathfrak{b}$ be a $(k+1)$-dimensional subspace, and let $\mathfrak{D}_k(\Lambda)$ be the degeneracy locus constructed in Section 3.6. The main properties of $\mathfrak{D}_k(\Lambda)$ are summarized in Corollary 13.

Recall that an irreducible component of $\mathfrak{D}_k(\Lambda)$ is said to be *standard* if it is generically supported on $X_0$. All the other irreducible components are *exceptional*.

LEMMA 20. *For a sufficiently general $\Lambda$ and a positive $k$, there is exactly one standard component of $\mathfrak{D}_k(\Lambda)$, and this component is generically reduced.*

*Proof.* Over the open subset $X^\circ$, the logarithmic tangent bundle agrees with the usual tangent bundle. Therefore

$$\mathfrak{X}_{X^\circ} = \Sigma_{X^\circ}.$$

First we show that $\mathfrak{D}_k(\Lambda) \cap X_0$ is irreducible. Since $X^\circ$ has a point fixed by a maximal torus of $B$, Lemma 16 says that

$$\mathrm{pr}_{2,X^\circ} : \Sigma_{X^\circ} \longrightarrow \mathbb{P}(\mathfrak{b})$$

is a dominant morphism. Therefore Bertini's theorem applies to $\mathrm{pr}_{2,X^\circ}$ and positive dimensional linear subspaces of $\mathbb{P}(\mathfrak{b})$ [Laz04, Theorem 3.3.1]. It follows that

$$\mathfrak{D}_k(\Lambda) \cap X^\circ = \mathrm{pr}_{1,X^\circ}\left(\mathrm{pr}_{2,X^\circ}^{-1}\left(\mathbb{P}(\Lambda)\right)\right)$$

is irreducible for a sufficiently general $\Lambda$.

Next we show that $\mathfrak{D}_k(\Lambda) \cap X_0$ is reduced. The tangent bundle of $X^\circ$ is generated by global sections from $\mathfrak{b}$, and hence there is a morphism to the Grassmannian

$$\Psi : X^\circ \longrightarrow \mathrm{Gr}_d(\mathfrak{b}), \qquad x \longmapsto \mathfrak{b}_x \quad \text{where} \quad d = \dim B - \dim X.$$

As a scheme, $\mathfrak{D}_k(\Lambda) \cap X^\circ$ is the pull-back of the Schubert variety in $\mathrm{Gr}_d(\mathfrak{b})$ defined by $\Lambda$. Therefore $\mathfrak{D}_k(\Lambda) \cap X^\circ$ is reduced for a sufficiently general $\Lambda$ by Kleiman's transversality theorem [Kle74, Remark 7]. $\square$

In fact, $\mathfrak{D}_k(\Lambda)$ has no embedded components for a sufficiently general $\Lambda$ (being a degeneracy locus of the expected dimension $k$), but we will not need this. When $Y$ is the Schubert variety $\mathbb{S}(\underline{\alpha})$, the reduced image in $\mathbb{S}(\underline{\alpha})$ of the unique standard component of $\mathfrak{D}_k(\Lambda)$ will be the subvariety $Z(\underline{\alpha})$ of Theorem 2.

*Proof of Theorem 15.* When $k$ is positive, there is exactly one standard component by Lemma 20. Write $\pi_*$ for the push-forward

$$\pi_* : A_*(X) \longrightarrow A_*(Y).$$

Our goal is to show that $\pi_*[\mathfrak{E}] = 0$ for all exceptional components $\mathfrak{E}$ of $\mathfrak{D}_k(\Lambda)$, for a sufficiently general $\Lambda$.

For this we consider the case when $k = 0$. Recall from Corollary 13 that $\mathfrak{D}_0(\Lambda)$ consists of finite set of points, each contained in a $B$-orbit $S$ such that $\mathfrak{X}_S = \Sigma_S$, for a sufficiently general

$\Lambda$. The number of points in $\mathfrak{D}_0(\Lambda)$ is equal to

$$\chi(X^\circ) = \int_X c_{SM}(X^\circ) = \sum_S \deg\left(\mathrm{pr}_{2,S} : \Sigma_S \longrightarrow \mathbb{P}(\mathfrak{b})\right) = 1,$$

where the sum is over all orbits such that $\mathfrak{X}_S = \Sigma_S$. Together with Lemma 16, the formula shows that every such orbit, except one, is of the form

$$S \simeq B/H, \qquad \mathrm{rank}(B) > \mathrm{rank}(H).$$

This one exception should be $X^\circ$, because $X^\circ$ contains a point fixed by a maximal torus of $B$.

Return to the case when $k$ is positive. Let $S$ be an orbit with $\mathfrak{X}_S = \Sigma_S$, and suppose that $S$ is different from $X^\circ$. Consider the $B$-equivariant map

$$\pi|_S : S \simeq B/H \longrightarrow \pi(S), \qquad \mathrm{rank}(B) > \mathrm{rank}(H).$$

The image of $S$ contains a point fixed by a maximal torus of $B$, because it is a $B$-orbit in $G/P$. Therefore $\pi(S)$ is of the form

$$\pi(S) \simeq B/K, \qquad \mathrm{rank}(B) = \mathrm{rank}(K).$$

Since $\pi$ is a regular log-resolution, this shows that Lemma 19 applies to $\pi|_S$. The degeneracy locus $D_k(\Lambda)$ of Lemma 19 is precisely the intersection $S \cap \mathfrak{D}_k(\Lambda)$ in our case because $\mathfrak{X}_S = \Sigma_S$. The conclusion is that

$$\dim \mathfrak{E} > \dim \pi(\mathfrak{E})$$

for any irreducible component $\mathfrak{E}$ of $\mathfrak{D}_k(\Lambda)$ generically supported on $S$.

Therefore $\pi_*[\mathfrak{E}] = 0$ for all exceptional components $\mathfrak{E}$, for a sufficiently general $\Lambda$. $\square$

## 5. A regular resolution of a classical Schubert variety

In this section, $E$ is a vector space with an ordered basis $e_1, \ldots, e_{n+d}$, $G$ is the general linear group of $E$, and $B$ is the subgroup of $G$ which consists of all invertible upper triangular matrices with respect to the ordered basis of $E$.

### 5.1

We recall the known resolution of singularities of the classical Schubert variety $\mathbb{S}(\underline{\alpha})$ which is regular in the sense of Definition 14. Theorem 2 therefore can be deduced from Theorem 15.

Let $\underline{\alpha} = (\alpha_1 \geqslant \alpha_2 \geqslant \cdots \geqslant \alpha_d \geqslant 0)$, and let $\mathbb{S}(\underline{\alpha}) \subseteq \mathrm{Gr}_d(E)$ be the Schubert variety defined with respect to the complete flag

$$F_\bullet = \left(F_0 \subsetneq F_1 \subsetneq \cdots \subsetneq F_{n+d}\right) \quad \text{where} \quad F_k := \mathrm{span}(e_1, \ldots, e_k).$$

DEFINITION 21. $\mathbb{V}(\underline{\alpha})$ is the subvariety

$$\mathbb{V}(\underline{\alpha}) := \left\{V_1 \subsetneq V_2 \subsetneq \cdots \subsetneq V_d \mid V_i \subseteq F_{\alpha_{d+1-i}+i}\right\} \subseteq \mathrm{Gr}_1(E) \times \mathrm{Gr}_2(E) \times \cdots \times \mathrm{Gr}_d(E).$$

The restriction to $\mathbb{V}(\underline{\alpha})$ of the projection to $\mathrm{Gr}_d(E)$ will be written

$$\pi_{\underline{\alpha}} : \mathbb{V}(\underline{\alpha}) \longrightarrow \mathbb{S}(\underline{\alpha}).$$

The projection $\pi_{\underline{\alpha}}$ maps $\mathbb{V}(\underline{\alpha})$ into $\mathbb{S}(\underline{\alpha})$ because $V_i \subseteq V_d \cap F_{\alpha_{d+1-i}+i}$ for all $i$.

We note that $\pi_{\underline{\alpha}}$ is the resolution used in [KL74] to obtain the determinantal formula for the classes of Schubert schemes. This resolution was also used in [AM09] to compute the Chern-Schwartz-MacPherson class of $\mathbb{S}(\underline{\alpha})^\circ$. All the properties of $\pi_{\underline{\alpha}}$ we need can be found in [AM09, Section 2]. However, one simple but important point for us was not emphasized in the non-embedded description of $\mathbb{V}(\underline{\alpha})$ in [AM09] as a tower of projective bundles: $\mathbb{V}(\underline{\alpha})$ is a subvariety of the partial flag variety

$$\mathrm{Fl}_{1,\ldots,d}(E) \subseteq \mathrm{Gr}_1(E) \times \mathrm{Gr}_2(E) \times \cdots \times \mathrm{Gr}_d(E),$$

and $\mathbb{V}(\underline{\alpha})$ is invariant under the diagonal action of $B$. It follows that

(i) $\mathbb{V}(\underline{\alpha})$ has finitely many $B$-orbits, and
(ii) every $B$-orbit of $\mathbb{V}(\underline{\alpha})$ contains a point fixed by a maximal torus of $B$.

The above properties imply that $\pi_{\underline{\alpha}}$ is a regular log-resolution of $\mathbb{S}(\underline{\alpha})$ in the sense of Definition 14.

*Remark* 22. We note that the Bott-Samelson variety of [Dem74, Han73] will not have finitely many $B$-orbits in general. It would be interesting to know which Schubert varieties in flag varieties (do not) admit a regular or $B$-finite log-resolution.

**5.2**

For the sake of completeness, we give an argument here that $\pi_{\underline{\alpha}}$ is a regular log-resolution of singularities of $\mathbb{S}(\underline{\alpha})$.

PROPOSITION 23. *$\pi_{\underline{\alpha}}$ is a regular log-resolution of $\mathbb{S}(\underline{\alpha})$. That is,*

(i) *$\pi_{\underline{\alpha}}$ is proper and $B$-equivariant,*
(ii) *$\pi_{\underline{\alpha}}^{-1}(\mathbb{S}(\underline{\alpha})^\circ) \longrightarrow \mathbb{S}(\underline{\alpha})^\circ$ is an isomorphism,*
(iii) *$\mathbb{V}(\underline{\alpha})$ is smooth and has finitely many $B$-orbits,*
(iv) *the complement of $\pi_{\underline{\alpha}}^{-1}(\mathbb{S}(\underline{\alpha})^\circ)$ in $\mathbb{V}(\underline{\alpha})$ is a divisor with normal crossings, and*
(v) *the isotropy Lie algebra $\mathfrak{b}_x$ contains a regular element of $\mathfrak{b}$ for each $x \in \mathbb{V}(\underline{\alpha})$.*

*Proof.* We start by justifying (ii). Note that $\pi_{\underline{\alpha}}$ has a section over the Schubert cell

$$s_{\underline{\alpha}} : \mathbb{S}(\underline{\alpha})^\circ \longrightarrow \pi_{\underline{\alpha}}^{-1}(\mathbb{S}(\underline{\alpha})^\circ), \qquad V \longmapsto V \cap \left(F_{\alpha_d+1} \subsetneq F_{\alpha_{d-1}+2} \subsetneq \cdots \subsetneq F_{\alpha_1+d}\right).$$

The statement

$$s_{\underline{\alpha}} \circ \pi_{\underline{\alpha}}|_{\pi_{\underline{\alpha}}^{-1}(\mathbb{S}(\underline{\alpha})^\circ)} = \mathrm{id}_{\pi_{\underline{\alpha}}^{-1}(\mathbb{S}(\underline{\alpha})^\circ)}$$

is equivalent to the assertion that

$$V_i = V_d \cap F_{\alpha_{d+1-i}+i}$$

for all $i$ and for all $V_\bullet \in \mathbb{V}(\underline{\alpha})$ with $V_d \in \mathbb{S}(\underline{\alpha})$. This is clear because $V_i$ is contained in the right-hand side and the dimensions of both sides are the same. Therefore

$$\pi_{\underline{\alpha}}^{-1}(\mathbb{S}(\underline{\alpha})^\circ) \longrightarrow \mathbb{S}(\underline{\alpha})^\circ$$

is an isomorphism, proving (ii).

We prove (iii) by induction on the number of entries of $\underline{\alpha}$. Define

$$\widetilde{\underline{\alpha}} := (\alpha_2 \geqslant \alpha_3 \geqslant \cdots \geqslant \alpha_d \geqslant 0)$$

and consider the corresponding subvariety

$$\mathbb{V}(\widetilde{\underline{\alpha}}) \subseteq \mathrm{Gr}_1(E) \times \mathrm{Gr}_2(E) \times \cdots \times \mathrm{Gr}_{d-1}(E).$$

Restricting the projection map which forgets the last coordinate, we have

$$\mathrm{pr}_{\hat{d}} : \mathbb{V}(\underline{\alpha}) \longrightarrow \mathbb{V}(\widetilde{\underline{\alpha}}).$$

Let $\mathscr{F}_\bullet$ be the flag of trivial vector bundles over $\mathbb{V}(\widetilde{\underline{\alpha}})$ modeled on the flag of subspaces $F_\bullet$. Then we may identify $\mathrm{pr}_{\hat{d}}$ with the projective bundle

$$\mathbb{P}(\mathscr{F}_{\alpha_1+d}/\mathscr{V}_{d-1}) \longrightarrow \mathbb{V}(\widetilde{\underline{\alpha}}),$$

where $\mathscr{V}_{d-1}$ is the pull-back of the tautological bundle from the projection $\mathbb{V}(\widetilde{\underline{\alpha}}) \longrightarrow \mathrm{Gr}_{d-1}(E)$. This shows by induction that $\mathbb{V}(\underline{\alpha})$ is smooth. The fact that $\mathbb{V}(\underline{\alpha})$ has finitely many $B$-orbits is implied by the Bruhat decomposition of $G$.

(iv) can also be proved by the same induction. Let $\widetilde{\underline{\alpha}}$ be as above, and set

$$D_{\mathrm{old}} := \mathbb{V}(\widetilde{\underline{\alpha}}) \setminus \pi_{\widetilde{\underline{\alpha}}}^{-1}\big(\mathbb{S}(\widetilde{\underline{\alpha}})^\circ\big).$$

We may suppose that $D_{\mathrm{old}}$ is a divisor in $\mathbb{V}(\widetilde{\underline{\alpha}})$ with normal crossings. The key observation is that

$$\mathbb{V}(\underline{\alpha}) \setminus \pi_{\underline{\alpha}}^{-1}\big(\mathbb{S}(\underline{\alpha})^\circ\big) = \mathrm{pr}_{\hat{d}}^{-1}(D_{\mathrm{old}}) \cup D_{\mathrm{new}},$$

where $D_{\mathrm{new}}$ is the smooth and irreducible divisor

$$D_{\mathrm{new}} := \mathbb{P}(\mathscr{F}_{\alpha_1+d-1}/\mathscr{V}_{d-1}) \subseteq \mathbb{P}(\mathscr{F}_{\alpha_1+d}/\mathscr{V}_{d-1}) = \mathbb{V}(\underline{\alpha}).$$

The assertion that $\mathrm{pr}_{\hat{d}}^{-1}(D_{\mathrm{old}}) \cup D_{\mathrm{new}}$ has normal crossings can be checked locally. Covering $\mathbb{V}(\underline{\alpha})$ with open subsets of the form $\mathrm{pr}_{\hat{d}}^{-1}(U)$, where $U$ is an open subset of $\mathbb{V}(\widetilde{\underline{\alpha}})$ over which the vector bundle $\mathscr{V}_{d-1}$ is trivial, the assertion becomes clear.

(v) is a consequence of the fact that each $B$-orbit of $\mathbb{V}(\underline{\alpha})$ contains a point fixed by a maximal torus of $B$. It follows that every point of $\mathbb{V}(\underline{\alpha})$ is fixed by a maximal torus of $B$. Therefore all the isotropy Lie algebras contain a Cartan subalgebra of $\mathfrak{b}$, whose general member is a regular element of $\mathfrak{b}$. □

## References

Alu99    Paolo Aluffi, *Differential forms with logarithmic poles and Chern-Schwartz-MacPherson classes of singular varieties*, Comptes Rendus de l'Académie des Sciences. Série I. Mathématique **329** (1999), no. 7, 619–624.

Alu05    Paolo Aluffi, *Characteristic classes of singular varieties*, Topics in Cohomological Studies of Algebraic Varieties, 1–32, Trends in Mathematics, Birkhäuser, Basel, 2005.

Alu06a    Paolo Aluffi, *Classes de Chern des variétés singulières, revisitées*, Comptes Rendus Mathématique. Académie des Sciences. Paris **342** (2006), no. 6, 405–410.

Alu06b    Paolo Aluffi, *Limits of Chow groups, and a new construction of Chern-Schwartz-MacPherson classes*, Pure and Applied Mathematics Quaterly **2** (2006), no. 4, 915–941.

AM09    Paolo Aluffi and Leonardo Constantin Mihalcea, *Chern classes of Schubert cells and varieties*, Journal of Algebraic Geometry **18** (2009), no. 1, 63–100.

Bor91    Armand Borel, *Linear Algebraic Groups*, Second edition, Graduate Texts in Mathematics **126**, Springer-Verlag, New York, 1991.

Bry10    Robert Bryant, *Rigidity and Quasi-Rigidity of Extremal Cycles in Hermitian Symmetric Spaces*, Annals of Mathematics Studies **153**, Princeton University Press, 2010.

BSS09    Jean-Paul Brasselet, José Seade, and Tatsuo Suwa, *Vector Fields on Singular Varieties*, Lecture Notes in Mathematics **1987**, Springer-Verlag, Berlin-Heidelberg, 2009.

Bri05    Michel Brion, *Lectures on the geometry of flag varieties*, Topics in Cohomological Studies of Algebraic Varieties, 33–85, Trends in Mathematics, Birkhäuser, Basel, 2005.

Bri07    Michel Brion, *Log homogeneous varieties*, Proceedings of the XVIth Latin American Algebra Colloquium, 1–39, Biblioteca de la Revista Matemática Iberoamericana Revista Matemática Iberoamericana, Madrid, 2007.

Bri09    Michel Brion, *Vanishing theorems for Dolbeault cohomology of log homogeneous varieties*, Tohoku Mathematical Journal (2) **61** (2009), no. 3, 365–392.

BJ08    Michel Brion and Roy Joshua, *Equivariant Chow ring and Chern classes of wonderful symmetric varieties of minimal rank*, Transformation Groups **13** (2008), 471–493.

BK05    Michel Brion and Ivan Kausz, *Vanishing of top equivariant Chern classes of regular embeddings*, Asian Journal of Mathematics **9** (2005), no. 3-4, 489–496.

Cos11    Izzet Coskun, *Rigid and non-smoothable Schubert classes*, Journal of Differential Geometry **87** (2011), no. 3, 493–514.

CR13    Izzet Coskun and Colleen Robles, *Flexibility of Schubert Classes*, preprint, arXiv:1303.0253.

Del70    Pierre Deligne, *Équations différentielles à points singuliers réguliers*, Lecture Notes in Mathematics **163**, Springer-Verlag, Berlin-New York, 1970.

Dem74    Michel Demazure, *Désingularisation des variétés de Schubert généralisées*, Annales Scientifiques de l'École Normale Supérieure (4) **7** (1974), 53–88.

FMSS95    William Fulton, Robert MacPherson, Frank Sottile, and Bernd Sturmfels, *Intersection theory on spherical varieties*, Journal of Algebraic Geometry **4** (1995), no. 1, 181–193.

Ful97    William Fulton, *Young Tableaux*, London Mathematical Society Student Texts **35**, Cambridge University Press, 1997.

Ful98    William Fulton, *Intersection Theory*, Second Edition. Ergebnisse der Mathematik und ihrer Grenzgebiete. 3. Folge. A Series of Modern Surveys in Mathematics **2**, Springer-Verlag, Berlin, 1998.

GP02    Mark Goresky and William Pardon, *Chern classes of automorphic vector bundles*, Inventiones Mathematicae **147** (2002), 561–612.

Han73    H. C. Hansen, *On cycles in flag manifolds*, Mathematica Scandinavica **33** (1973), 269–274.

Hon05    Jaehyun Hong, *Rigidity of singular Schubert varieties in $Gr(m,n)$*, Journal of Differential Geometry **71** (2005), no. 1, 1–22.

Hon07    Jaehyun Hong, *Rigidity of smooth Schubert varieties in Hermitian symmetric spaces*, Transactions of the American Mathematical Society **359** (2007), no. 5, 2361–2381.

Huh12a    June Huh, *Milnor numbers of projective hypersurfaces and the chromatic polynomial of graphs*, Journal of the American Mathematical Society **25** (2012), 907–927.

Huh12b    June Huh, *h-vectors of matroids and logarithmic concavity*, preprint, `arXiv:1201.2915`.

Jon07    Benjamin Jones, *On the singular Chern classes of Schubert varieties via small resolution* Ph. D. Thesis, University of Notre Dame, 2007.

Jon10    Benjamin Jones, *Singular Chern classes of Schubert varieties via small resolution*, International Mathematics Research Notices (2010), no. 8, 1371–1416.

KL74    George Kempf and Dan Laksov, *The determinantal formula of Schubert calculus*, Acta Mathematica **132** (1974), 153–162.

Ken90    Gary Kennedy, *MacPherson's Chern classes of singular algebraic varieties*, Communications in Algebra **18** (1990), no. 9, 2821–2839.

Kir06    Valentina Kiritchenko, *Chern classes of reductive groups and an adjunction formula*, Annales de l'Institut Fourier (Grenoble) **56** (2006), no. 4, 1225–1256.

Kir07    Valentina Kiritchenko, *On intersection indices of subvarieties in reductive groups*, Moscow Mathematical Journal **7** (2007), no. 3, 489–505.

Kle74    Steven Kleiman, *The transversality of a general translate*, Compositio Mathematica **28** (1974), 287–297.

Laz04    Robert Lazarsfeld, *Positivity in Algebraic Geometry I*, Ergebnisse der Mathematik und ihrer Grenzgebiete **48**, Springer-Verlag, Berlin, 2004.

Mac74    Robert MacPherson, *Chern classes for singular algebraic varieties*, Annals of Mathematics (2) **100** (1974), 423–432.

MO67    Hideyuki Matsumura and Frans Oort, *Representability of group functors, and automorphisms of algebraic schemes*, Inventiones Mathematicae **4** (1967), 1–25.

Mih07    Leonardo Constantin Mihalcea, *Sums of binomial determinants, non-intersecting lattice paths, and positivity of Chern-Schwartz-MacPherson classes*, preprint, arXiv:0702566.

Per02    Nicolas Perrin, *Courbes rationnelles sur les variétés homogènes*, Annales de l'Institut Fourier (Grenoble) **52** (2002), no.1, 105–132.

Ram64    C. P. Ramanujam, *A note on automorphism groups of algebraic varieties*, Mathematische Annalen **156** (1964), 25-33.

Sai80    Kyoji Saito, *Theory of logarithmic differential forms and logarithmic vector fields*, Journal of the Faculty of Science. University of Tokyo. Section IA. Mathematics **27** (1980), no. 2, 265–291.

Sch05    Jörg Schürmann, *Lectures on characteristic classes of constructible functions*, Notes by Piotr Pragacz and Andrzej Weber, Trends in Mathematics, Topics in Cohomological Studies of Algebraic Varieties, 175–201, Birkhäuser, Basel, 2005.

Sch65a    Marie-Hélène Schwartz, *Classes caractéristiques définies par une stratification d'une variété analytique complexe. I.* Comptes Rendus Mathématique. Académie des Sciences. Paris **260** (1965), 3262–3264.

Sch65b    Marie-Hélène Schwartz, *Classes caractéristiques définies par une stratification d'une variété analytique complexe. II.* Comptes Rendus Mathématique. Académie des Sciences. Paris **260** (1965) 3535–3537.

Str11    Judson Stryker, *Chern-Schwartz-MacPherson classes of graph hypersurfaces and Schubert varieties*, Ph. D. Thesis, Florida State University, 2011.

TY05    Patrice Tauvel and Rupert Yu, *Lie Algebras and Algebraic Groups*, Springer Monographs in Mathematics, Springer-Verlag, Berlin, 2005.

Web12    Andrzej Weber, *Equivariant Chern classes and localization theorem*, Journal of Singularities **5** (2012), 153–176.

June Huh    junehuh@umich.edu

Department of Mathematics, University of Michigan, Ann Arbor, MI 48109, USA